\documentstyle{article}

\newtheorem{theorem}{Theorem}[section]    
\newtheorem{lemma}[theorem]{Lemma}         
\newtheorem{corollary}[theorem]{Corollary}

\def\PSL{\mbox{\rm{PSL}}}

\def\SL{\mbox{\rm{SL}}}
\def\SO{\mbox{\rm{SO}}}
\def\Isom{\mbox{Isom}}

\def\GL{\mbox{\rm{GL}}}

\def\qed{ $\Box$}
\def\O{\mbox{\rm{O}}}

\def\det{\mbox{\rm{det}}}
\def\stab{\mbox{\rm{Stab}}}

\def\mod{\mbox{mod}}
\def\demo{ {\bf Proof.} }

\title{The Bianchi groups are separable on geometrically finite subgroups.}
\author{I. Agol,~ D. D. Long\thanks{This work was partially supported by
the N. S. F. }  ~ \& A. W. Reid\thanks{This work was partially supported by 
the N. S. F., The Alfred P. Sloan Foundation and a grant
from the Texas Advanced Research Program.}}

\begin{document}           

\maketitle
\begin{abstract}
Let $d$ be a square free positive integer and $O_d$ the ring of integers
in ${\bf Q}(\sqrt{-d})$. The main result of this
paper is to show that the groups $\PSL(2,O_d)$ are
subgroup separable on geometrically finite subgroups.
\end{abstract}
\section{Introduction}

Let $G$ be a group and $H$ a finitely generated subgroup, $G$ is
called {\em $H$-subgroup separable} if given any $g \in G\setminus H$,
there exists a subgroup $K < G$ of finite index with $H < K$ and
$g\notin K$. $G$ is called {\em subgroup separable} (or {\em LERF}) if
$G$ is $H$-subgroup separable for all finitely generated $H < G$. 
Subgroup separability is an extremely powerful property, for
instance it is much stronger than residual finiteness. 
The class of groups for which subgroup separability is known for
{\em all} finitely generated subgroups is extremely small;
abelian groups, free groups, surface groups and carefully controlled 
amalgamations of these, see \cite{Gi}, and \cite{Sc} for example.    

However our motivation comes from 3-manifold topology, where the importance of 
subgroup separability stems from the well-known fact
(cf. \cite{Sc} and \cite{Lo1}) that it allows passage from immersed 
incompressible surfaces to embedded incompressible surfaces in finite covers.
It therefore makes sense (especially in light of the facts that there are 
closed 3-manifolds $M$ for which $\pi_1(M)$ is {\em not}
subgroup separable \cite{BKS}) to ask for separability only for some
mildly restricted class of subgroups.
  
No example of a finite co-volume
Kleinian group is known to be subgroup separable. However, in this
context the geometrically finite subgroups (especially the
geometrically finite surface subgroups) are the intractable and most
relevant case in all applications.  The reason for this is that the work of Bonahon and
Thurston (See \cite{Bo}) implies that freely indecomposable
geometrically infinite subgroups of finite co-volume Kleinian groups
are virtual fiber groups, and these are easily seen to be separable.
Accordingly, there has been much more attention paid to separating
geometrically finite subgroups of finite co-volume Kleinian groups
(see \cite{Gi}, \cite{Wi}). A class of Kleinian groups that have been
historically important in the subject are the Bianchi groups.  Our
main result is the following:

\begin{theorem}
\label{main}
Let $d$ be a square free positive integer, and $O_d$ the
ring of integers in ${\bf Q}(\sqrt{-d})$. The Bianchi group
$\PSL(2,O_d)$ is $H$-subgroup separable for all geometrically finite subgroups
$H$.
\end{theorem}

A case which has attracted much interest itself is the fundamental group of
the figure eight knot complement. This has index $12$ in $\PSL(2,O_3)$. Hence
we get (see also D. Wise  \cite{Wi} in this case):

\begin{corollary}
\label{fig8}
Let $K$ denote the figure eight knot, then $\pi_1(S^3 \setminus K)$ is
$H$-subgroup separable for all geometrically finite subgroups $H$.
\qed\end{corollary}

The fundamental group of the
Borromean rings is well-known \cite{Th} to be a subgroup of index $24$ in
$\PSL(2,O_1)$, hence we also have:

\begin{corollary}
The fundamental group of the Borromean rings is subgroup separable
on its geometrically finite subgroups.
\qed\end{corollary}

In fact, since it is easy to show that there are infinitely many (2 component)
links in $S^3$ whose complements are arithmetic
we deduce,

\begin{corollary}
There are infinitely many hyperbolic links in $S^3$ for which the
fundamental group of the complement is separable on all geometrically finite
subgroups.\qed\end{corollary}

Our methods also apply to give new examples of cocompact Kleinian groups
which are separable on all geometrically finite subgroups.

\begin{theorem}
\label{cocompact}
There exist infinitely many commensurability classes of 
cocompact Kleinian groups which satisfy:
\begin{itemize}
\item they are $H$-subgroup separable for all geometrically finite subgroups
$H$, 
\item they are not commensurable with a group generated by reflections. 
\end{itemize}
\end{theorem}

The second statement of Theorem \ref{cocompact} is only included to distinguish
the groups constructed from those that Scott's argument (\cite{Sc} and below)
applies to in dimension $3$.

In the cocompact setting some interesting groups that we can handle are the
following. Let $\Gamma$
denote the subgroup of index 2 consisting of orientation-preserving
isometries in the group generated by reflections in the faces of the
tetrahedron in ${\bf H}^3$ described as $T_2[2,2,3;2,5,3]$ (see \cite{MR} for
instance for notation).  This tetrahedron
has a symmetry of order $2$, and this symmetry extends to an
orientation preserving isometry of the orbifold $Q = {\bf H}^3/\Gamma$. Let
$\Gamma_0$ be the Kleinian group obtained as the orbifold group of this
2-fold quotient. $\Gamma_0$ attains the minimal co-volume
for an arithmetic Kleinian group (see below), \cite{CF}. It is conjectured
to attain the smallest co-volume for all Kleinian groups.  We show as 
corollaries of the methods,

\begin{corollary}
\label{cox}
$\Gamma_0$ is $H$-subgroup separable for all geometrically finite subgroups
$H$.\end{corollary}

\begin{corollary}
\label{SWD}
Let $W$ denote the Seifert-Weber dodecahedral space, then $\pi_1(W)$ is
$H$-subgroup separable for all geometrically finite subgroups $H$.
\end{corollary}

Although $\Gamma_0$ and $\pi_1(W)$ are commensurable with groups
generated by reflections, as far as we know they are not commensurable
with one where all dihedral angles are $\pi/2$, as is required in
applying \cite{Sc}.  Note that by \cite{Th} Chapter 13, the group
$\Gamma$ does not split as a free product with amalgamation or
HNN-extension since the orbifold ${\bf H}^3/\Gamma$ is non-Haken in
the language of orbifolds.  It is also widely believed that $W$ is
non-Haken, and these appear to the first explicit examples of such
Kleinian groups separable on geometrically finite subgroups (see also \S 6.1
for a further example).

An application of Theorem \ref{main} that
seems worth recording is the following.
An obvious subgroup of $\PSL(2,O_d)$ is $\PSL(2,{\bf Z})$, and in the
context of the Congruence Kernel, Lubotzky \cite{Lu} asked the
following question: \\[\baselineskip] 
{\bf Question}: Is the induced map:
 $$\eta_d~:~\widehat{\PSL}(2,{\bf Z}) \rightarrow
\widehat{\PSL}(2,O_d)$$ 
injective?   \\[\baselineskip]
Here $\hat{G}$ denotes the profinite completion (see \S 8 for details). By
a standard reformulation of $H$-subgroup separable, we are able to give
an affirmative answer.

\begin{theorem}
\label{inject}
The map $\eta_d$ is injective for all $d$.
\end{theorem} 
As is pointed out in \cite{Lu},
this has important ramifications
for the nature of the Congruence Kernel and
the structure of non-congruence subgroups of the Bianchi groups. For example
Theorem \ref{inject} gives another proof that the Bianchi groups do
not have the Congruence Subgroup Property \cite{Se}.

Here is an overview of the paper. The underlying method in proving 
Theorems \ref{main} and \ref{cocompact} 
is to use the arithmetic theory of quadratic forms and 
their relationship to discrete groups of hyperbolic isometries to inject (up
to commensurability) Kleinian groups into a fixed finite co-volume
group acting on a higher dimensional hyperbolic space 
commensurable with a group
generated by reflections in an all right polyhedron in ${\bf H}^n$.
This part of the proof hinges upon Lemmas \ref{subgpbianchi} and
\ref{infinitelymany}.

The ambient reflection group is then shown to be subgroup separable on its 
geometrically finite subgroups. This has its origins in \cite{Sc} where 
it is shown that in dimension $2$, cocompact Fuchsian groups are subgroup 
separable. In \S 3, we generalise this to all dimensions. The 
situation is actually a good deal more delicate than is generally 
appreciated; in
particular, in Scott's article \cite{Sc} various statements are made
that have been taken to suggest that his methods extend to higher
dimensions and that they could be used to separate geometrically
finite subgroups inside  Kleinian groups commensurable with groups 
generated by reflections in all right ideal polyhedra in ${\bf H}^3$. 
However it seems to us (and has been confirmed by Scott) that this is
not so.\\[\baselineskip]  
{\bf Acknowledgement:}~ The first author would like to thank 
Bill Thurston for some
useful conversations in regards to the proof of Lemma \ref{subgpbianchi}.

\section{Preliminaries}
\label{preliminaries}

We recall some facts about discrete groups of isometries of hyperbolic
spaces, and arithmetic hyperbolic n-manifolds, \cite{R} and
\cite{B}. A reasonable reference that contains information on both
these topics is \cite{VS}.

\subsection{}

Let $f_n$ be the $(n+1)$-dimensional quadratic form $<1,1\ldots ,1,-1>$.
The orthogonal group of this form is simply $O(n,1;{\bf R})$. This has
four connected components. Passing to $\SO(n,1;{\bf R})$, this has
two connected components, 
and the connected component of the identity in 
$\SO(n,1;{\bf R})$, denoted $\SO_0(f_n;{\bf R})$ (which has finite
index in $\O(n,1;{\bf R})$), may be identified with $\Isom_+({\bf H}^n)$;
it preserves the upper sheet of the hyperboloid $f_n(x) = -1$
and the orientation. 
Given a (discrete) subgroup $G$ of $\O(n,1;{\bf R})$, 
$G \cap \SO_0(n,1;{\bf R})$ has finite index in $G$. 

\subsection{}

A Kleinian group will always refer to a discrete group of orientation
preserving isometries of ${\bf H}^3$.  Throughout the paper we often pass
between models of ${\bf H}^3$, and use the term Kleinian group in both
settings. Hopefully no confusion will arise.

A Kleinian group $\Gamma$ is {\em geometrically finite} if either of the
following equivalent statements hold (see \cite{R}):
\begin{enumerate}
\item  $\Gamma$ admits a finite sided Dirichlet polyhedron for its action
on ${\bf H}^3$.
\item Let $C(\Gamma)$ denote the convex core of ${\bf H}^3/\Gamma$,
then for all $\epsilon > 0$ the $\epsilon$-neighbourhood,
$N_\epsilon(C(\Gamma))$ has finite volume.
\end{enumerate}
\smallskip
In higher dimensions, geometrical finiteness has been more problematical,
see \cite{Bw} for example. For instance, the generalization of 1 above
(which is the classical statement in dimensions $2$ and $3$) becomes
more restrictive in higher dimensions, cf. \cite{Bw}.
However, if we insist that $\Gamma$ is {\em finitely generated}, then the
statement 2 above suffices as a definition for $\Gamma$ to be geometrically
finite, see \cite{Bw}. Henceforth, the term geometrically finite always
refers to this. That finite generation is important follows from
work of \cite{Ha}. 

We will also make use of the following equivalent statement of geometrical
finiteness ($GF2$ of \cite{Bw}). This requires some terminology from
\cite{Bw}.

Suppose that $\Gamma$ is a discrete subgroup of $\Isom_+({\bf H}^n)$, and
$\Lambda(\Gamma)$ is its limit set.  Let $p\in \Lambda(\Gamma)$ be a
parabolic fixed point and 
$$\stab_p(\Gamma) = \{\gamma\in \Gamma~:~\gamma p = p\}.$$
$p$ (as above) is called {\em bounded} if 
$(\Lambda(\Gamma) \setminus p)/\stab_p(\Gamma)$ is compact. Then \cite{Bw}
shows:

\begin{lemma}
\label{GF2}
Suppose that $\Gamma$ is a discrete subgroup of $\Isom_+({\bf H}^n)$. 
Then $\Gamma$ is geometrically finite if and only if the limit
set consists of conical limit points and bounded parabolic fixed points.
\qed\end{lemma}

As is pointed out in \cite{Bw} the notion of bounded parabolic fixed point
can be stated as: $p$ is a bounded parabolic fixed point if and only if
given a minimal plane $\tau$ 
which is invariant for the action of $\stab_\Gamma(p)$ 
then there is a
constant $r$ with the property that $\Lambda(\Gamma) \setminus \{p\}$ 
lies inside an $r$-neighbourhood 
(measured in a Euclidean metric on ${\bf R}^{n-1}$) of the plane $\tau$.

\subsection{}

We will require some notions from the theory of quadratic forms, we
refer to  \cite{La} as a standard reference.

We first fix some notation, if $k$ is a number field, then $R_k$ will
denote its ring of integers.\\[\baselineskip]
{\bf Definition.} Two $n$-dimensional quadratic forms $f$
and $q$ defined over a field $k$ (with associated symmetric matrices
$F$ and $Q$) are {\em equivalent} over $k$ if there exists $P \in
\GL(n, k)$ with $P^tFP = Q$.\\[\baselineskip]
Let $f$ be a quadratic form in
$n+1$ variables with coefficients in a real number field $k$, with
associated symmetric matrix $F$, and let
$$\SO(f) = \{X \in \SL(n+1,{\bf C})~|~X^tFX = F \}$$ 
be the Special Orthogonal group of $f$.  This is
an algebraic group defined over $k$, and $\SO(f;R_k)$
is an arithmetic subgroup of $\SO(f;{\bf R})$, \cite{BH} or \cite{B}. In
particular $\SO(f;{\bf R})/\SO(f;R_k)$ has finite volume with respect
to a suitable measure.

\begin{lemma}
\label{equivalent}
In the notation above, 
$\SO(f;{\bf R})$ is conjugate to $\SO(q;{\bf R})$,
$\SO(f;k)$ is conjugate to $\SO(q;k)$, and
$\SO(f;R_k)$ is conjugate to a subgroup of $\SO(q;k)$
commensurable with $\SO(q;R_k)$.
\end{lemma}
\demo  We do the $k$-points case first, the real case follows
the same line of argument. Thus, let $A \in \SO(q;k)$, then,
$$(PAP^{-1})^tF(PAP^{-1}) = (P^{-1})^tA^t(P^tFP)AP^{-1} = 
(P^{-1})^tA^tQAP^{-1} .$$
Since $A \in \SO(q;k)$, this gives,
$$(P^{-1})^tQP^{-1} = F$$
by definition of $P$. Since $P$ has $k$-rational entries, we have that
$PAP^{-1} \in \SO(f;k)$. The reverse inclusion is handled
similarly.

To handle the second statement, observe that by considering the
denominators of the entries of $P$ and $P^{-1}$, we can choose a deep enough
congruence subgroup $\Gamma$ in $\SO(q;R_k)$ so that that $P\Gamma
P^{-1}$ will be a subgroup of $\SO(f;R_k)$. The index is
necessarily finite since both are of finite
co-volume.\qed\\[\baselineskip] 
Assume now that $k$ is totally real, and let $f$ be a form in $n+1$-variables
with coefficients in $k$, and 
be equivalent over $\bf R$ to the
form $f_n$ (as in \S 2.1).
Furthermore, if
$\sigma~:~k\rightarrow {\bf R}$ is a field embedding, then the form
$f^\sigma$ obtained by applying $\sigma$ to $f$ is defined over the real
number field $\sigma (k)$. We insist that for embeddings $\sigma \neq id$,
$f^\sigma$ is equivalent over $\bf R$ to the form
in $(n+1)$-dimensions, of signature $(n+1,0)$.
Since $f$ is equivalent over $\bf R$ to $f_n$,
it from follows Lemma \ref{equivalent} that
$\SO(f;{\bf R})$ is conjugate in $\GL(n+1,{\bf R})$ to $\SO(f_n;{\bf
R})$. From the discussion in \S 2.1, we deduce from \cite{BH} (or \cite{B})
that $\SO_0(f;R_k)$ defines an arithmetic lattice in $\Isom_+({\bf H}^n)$.
For $n$ odd, this gives only a subclass of 
arithmetic lattices in $\Isom_+({\bf H}^n)$ (see \cite{VS} pp 221--222).  

The group $\SO_0(f;R_k)$ (and hence the conjugate in $\Isom_+({\bf
H}^n)$) is cocompact if and only if the form $f$ does not represent
$0$ non-trivially with values in $k$, see \cite{BH}.  Whenever $n
\geq 4$, the lattices constructed above are non-cocompact if and
only if the form has rational coefficients, since it is well known
every quadratic
form over $\bf Q$ in at least $5$ variables represents $0$
non-trivially, see \cite{La}. 

We make some comments on arithmetic hyperbolic 3-manifolds
constructed in the above manner. It is a consequence of the results in
\cite{MR} that when $n = 3$, the above class of arithmetic Kleinian groups
coincide precisely with those arithmetic Kleinian groups containing
a non-elementary Fuchsian subgroup and for which the invariant trace field
is quadratic imaginary. As we will require it we state a corollary of the main
result of \cite{MR}.  This ties up the $\PSL(2)$ and $\SO(3,1)$ descriptions
of certain arithmetic Kleinian groups we will need.

\begin{theorem}
\label{macr}
Let $a$, $b$ and $c$ be integers with $a<0$ and $b, c > 0$. Let $q$ be the
quadratic form $<1,a,b,c>$. Then $\SO_0(q,{\bf Z})$ defines an arithmetic
Kleinian subgroup of $\PSL(2,{\bf C})$ with invariant trace-field 
${\bf Q}(\sqrt{abc})$ and invariant quaternion algebra with Hilbert
Symbol $\biggl({{-ac,-bc} \over {\bf Q}(\sqrt{abc})}\biggr)$.\qed\end{theorem}

\noindent Indeed the correspondence above is a bijective
correspondence between commensurability classes in the two models. 

\section{All right reflection groups and separability}
\label{rightreflections}

Suppose that $P$ is a compact or ideal polyhedron (ie at least one vertex
lies at the sphere-at-infinity) in ${\bf H}^n$ all
of whose dihedral angles are $\pi/2$.  Henceforth we call this
an {\em all right polyhedron}. Then the Poincar\'e polyhedron
theorem implies that the group generated by reflections in the co-dimension
one faces
of $P$ is discrete and a fundamental domain for its action is the
polyhedron $P$, that is to say, we obtain a tiling of hyperbolic
$n$-space by tiles all isometric to $P$.  Let the group so generated
be denoted by $G(P)$. 
\begin{theorem}
\label{allright}
The group $G(P)$ is H-subgroup separable for every finitely generated
geometrically finite subgroup $H < G(P)$.
\end{theorem}
It seems to be a folklore fact that this theorem follows easily from
the ideas contained in \cite{Sc}; however it seems to us that this is
not the case and we include a complete proof. One piece of terminology 
we require is the following. Any horospherical 
cusp cross-section of a hyperbolic n-orbifold of
finite volume is finitely covered by the $n-1$-torus (see \cite{R} Chapter 5). 
We say a cusp of 
a hyperbolic n-orbifold is of {\em full rank} if it contains 
${\bf Z}^{n-1}$ as a subgroup of finite index. 
Otherwise the cusp is said not to be of full rank.\\[\baselineskip]
{\bf Proof of 3.1}
The proof breaks up into various cases which we deal with in ascending
order of difficulty. All proofs hinge upon the observation (see
\cite{Sc}) that the separability of $H$ is equivalent to the
following:\\[\baselineskip]
Suppose that we are given a compact subset $X \subset {\bf H}^n/H$. Then
there is a finite index subgroup $K < G(P)$, with $H < K$ and with
the projection map $q : {\bf H}^n/H \longrightarrow {\bf H}^n/K$ being an
embedding on $X$.

We sum up the common strategy which achieves this. The group $H$ is
geometrically finite and one can enlarge its convex hull in ${\bf
H}^n/H$ so as to include the compact set $X$ in a convex set contained
in ${\bf H}^n/H$; this convex set lifts to an $H$-invariant convex set
inside ${\bf H}^n$. One then defines a coarser convex hull using only
the hyperbolic halfspaces bounded by totally geodesic planes which
come from the $P$-tiling of ${\bf H}^n$; this hull is denoted by
$H_P(C^+)$.  This hull is $H$-invariant and the key point is to show that
$H_P(C^+)/H$ only involves a finite number of tiles. It is the
mechanics of achieving this that vary depending on the nature of $P$
and $H$; the remainder of the proof follows \cite{Sc} and is an
elementary argument using the Poincar\'e polyhedron theorem and some
covering space theory. The details are included in the first argument
below.

\subsection{ $P$ is compact. }
Let $C$ be a very small neighbourhood of the convex hull of $H$,
regarded as a subset of ${\bf H}^n$. In our setting, the group
$G(P)$ contains no parabolic elements so that the hypothesis implies
that $C/H$ is compact.

The given set $X$ is compact so that there is a $t$ with the property
that every point of $X$ lies within a distance $t$ of $C/H$. Let $C^+$ be the $10t$
neighbourhood of $C$ in ${\bf H}^n$. This is still a convex
$H$-invariant set and $C^+/H$ is a compact convex set containing $X$.

As discussed above, take the convex hull $H_P(C^+)$
 of $C^+$ in ${\bf H}^n$ using the half
spaces coming from the $P$-tiling of ${\bf H}^n$.
By construction $H_P(C^+)$ is a union of $P$-tiles, is convex and
$H$-invariant.  The crucial claim is:\\[\baselineskip]
\leftline{{\bf Claim}:~$H_P(C^+)/H$ involves only
 a finite number of such tiles.}\\[\baselineskip]
To see this we argue as follows. 

Fix once and for all a point in the interior of a top dimensional face
of the tile and call this its {\em barycentre}. The tiles we use actually
often have a geometric barycentre (i.e. a point which is equidistant
from all of the faces) but such special geometric properties are not used;
it is just a convenient reference point.

Our initial claim is that if the barycentre of
a tile is too far away from $C^+$, then it cannot lie in $H_P(C^+)$.

The reason for this is the convexity of $C^+$. If $a$ is a point in ${\bf H}^n$
not lying in  $C^+$ then  there is a unique point on $C^+$ which
is closest to $a$. Moreover, if this distance is $R$, then the set of points
distance precisely $R$ from $a$ is a sphere touching $C^+$ at a single point
$p$ on the frontier of $C^+$ and the geodesic hyperplane tangent to the sphere 
at this point is the (generically unique) supporting hyperplane separating 
$C^+$ from $a$. 

Suppose then that $P^*$ is a tile whose barycentre is very distant from
$C^+$.  Let $a^*$ be the point of $P^*$ which is closest to $C^+$ and
let $p$ be a point on the frontier of $C^+$ which is closest to $P^*$.
As noted above, there is a geodesic supporting hyperplane ${\cal H}_p$
through $p$ which is (generically) tangent to $C^+$ and separates
$C^+$ from $a^*$.  Let the geodesic joining $a^*$ and $p$ be denoted
by $\gamma$. Note that since $p$ is the point of $C^+$ closest to
$a^*$, $\gamma$ is orthogonal to ${\cal H}_p$.

If $a^*$ happens to be in the interior of a tile face of $P^*$, then this
tile face must be at right angles to $\gamma$, since $a^*$ was
closest. Let ${\cal H}_{a^*}$ be the tiling plane defined by this tile face.
Since in this case $\gamma$ is at right angles to both  ${\cal H}_{a^*}$
and ${\cal H}_{p}$, these planes are disjoint and so the tiling plane 
separates $P^*$ from $C^+$ as required.

If $a^*$ is in the interior of some smaller dimensional face,
$\sigma$, then the codimension one faces of $P^*$ which are incident
at $\sigma$ cannot all make small angles with $\gamma$ since they make
right angles with each other. The hyperplane  $\cal H$ which makes 
an angle close to $\pi/2$ plays the role of ${\cal H}_{a^*}$ 
in the previous paragraph. The reason is that since $a^*$ and $p$ are 
very distant and the planes ${\cal H}_p$ and $\cal H$ both make angles 
with $\gamma$ which are close to $\pi/2$, the planes are disjoint 
and we see as above that $P^*$ cannot lie in the tiling hull in this
case either.

The proof of the claim now follows, as there can be only 
finitely many barycentres near to any compact subset of ${\bf H}^n/H$.

 The proof of subgroup separability now finishes off as in
\cite{Sc}. Let $K_1$ be the subgroup of $G(P)$ generated by
reflections in the sides of $H_P(C^+)$.  The Poincar\'e polyhedron
theorem implies that $H_P(C^+)$ is a (noncompact) fundamental domain
for the action of the subgroup $K_1$. Set $K$ to be the subgroup of
$G(P)$ generated by $K_1$ and $H$, then ${\bf H}^n/K = H_P(C^+)/H$ so
that $K$ has finite index in $G(P)$. Moreover, the set $X$ embeds as
required.  $\Box$

\subsection{ $P$ is an ideal all right polyhedron.}
{\bf Subcase A:} $H$ has no cusps.\\[\baselineskip]
This case is very similar to the case that $G(P)$ is cocompact since
in the absence of cusps, the core of $H$ is actually compact. We form
the set $C^+$ as above.  The set $C^+/H$ is still compact so that it
only meets a finite number of tiles and we choose a constant $K$ so that
the barycentre of each such tile is within distance $K$ of $C^+/H$. 

Now we repeat the argument above, with the extra care that one should
only look at tiles in ${\bf H}^n$ whose barycentres are at distance from $C^+$
much larger than $K$; this ensures that such a tile cannot meet $C^+$  and
the rest of  the argument is now identical.\\[\baselineskip] 
{\bf Subcase B:} $H$ has cusps which are all of full rank.\\[\baselineskip]
In this case the core $C^+/H$ is no longer compact, but by geometrical
finiteness it has finite volume. The thick part of this core is
compact and can be covered by a finite number of tiles. Also the thin
part can be covered by a finite number of tiles; one sees this by
putting the cusp of $H$ at infinity, the cusp has full rank so there
is a compact fundamental domain for its action and this fundamental
domain meets only a finite number of tiles.

Choose $K$ for this finite collection of tiles, then argue as in
Subcase A.\\[\baselineskip] 
{\bf Subcase C:} $H$ has a cusp of less than full rank.
\\[\baselineskip]
The idea in this case is to enlarge $H$ to a group $H^*$ which now
only has full rank cusps in such a way that the compact set $X$
continues to embed in the quotient ${\bf H}^n/H^*$; we then argue as
in Subcase B. The argument follows a line established in \cite{CL}. We
assume that $H$ has a single cusp of less than full rank, the case of
many cusps is handled by successive applications of this case.

To this end, consider the upper $n$-space model for ${\bf H}^n$
arranged so that $\infty$ is a parabolic fixed point for $H$. Denote
the limit set of $H$ by $\Lambda(H)$.  By Lemma \ref{GF2} and the remarks
following it, if 
$\tau$ is a minimal plane which is invariant for the action of
$\stab_H(\infty)$ then there is a constant $r$ with the property that
$\Lambda(H) - \infty$ lies inside an $r$-neighbourhood (measured in a
Euclidean metric on ${\bf R}^{n-1}$) of the plane $\tau$.

We now sketch the construction of \cite{CL}. We have already observed
that $C^+/H$ is a finite volume hyperbolic $n$-orbifold; one can
therefore define a ``thickness" for this hull, that is to say, there
is a constant $c_1$ so that every point on the upper hypersurface of
$C^+$ is within distance at most $c_1$ of some point on the lower
hypersurface. We may as well suppose that $c_1$ is fairly large, say
at least $10$.

Choose a horoball $N'$ in ${\bf H}^n/H$ which is so small that it
is a very long way from the original compact set $X$. By shrinking
further, we arrange that the distance between any two preimages of
$N'$ in ${\bf H}^n$ is at least, say $1000c_1$.  Now shrink further
and find a horoball $N \subset N'$ so that $\partial N$ is distance
$1000c_1$ from $\partial N'$.  It follows that when we look in the
universal covering, if a preimage of $N$ is not actually centred on
some point of $\Lambda(H)$, then it is distance at least $750c_1$ from
$C^+$.

Using geometrical finiteness in the form given by Lemma \ref{GF2},
we may find  a pair of $\stab_H(\infty)$-invariant, parallel,
totally geodesic hyperbolic $n-1$-planes, both passing through $\infty$
which are distance $10r$ apart in the Euclidean metric on ${\bf R}^{n-1}$
and contain $\Lambda(H)- \infty$. By moving the planes further apart
if necessary, we may assume that the slab of hyperbolic $n$-space
between them contains all of $C^+$.  

If we denote the dimension of the minimal plane $\tau$ by $n-1-k$, then we 
may choose $k$ such pairs of planes so that the union of these of $2k$ planes 
cuts out subset of ${\bf R}^{n-1}$ which has the form $[0,1]^k \times {\bf R}^{n-1-k}$
 containing all of $\Lambda(H) - \infty$
and so that the slab of hyperbolic $n$-space cut out by these planes contains
all of $C^+$.  Denote this slab by $\Sigma$.

Let $N_\infty$ be the preimage of $N$ centred at $\infty$; $\partial
N_\infty$ is a copy of ${\bf R}^{n-1}$ equipped with a Euclidean metric
coming from the restriction of the metric coming from ${\bf H}^n$.  By
choosing translations in $\stab_{G(P)}(\infty)$ which move points a
very long distance in the Euclidean metric on $\partial N_\infty$, we
can augment $\stab_H(\infty)$ to form a new subgroup $\stab_H(\infty)^*
\subset \stab_{G(P)}(\infty)$, now of full rank, with the property that
any element of $\stab_H(\infty)^*$ either stabilises the slab $\Sigma$
(this is the case that the translation in question in fact lies in the
subgroup $\stab_H(\infty)$) or moves it a very long distance from
itself, measured in the Euclidean metric of $\partial N_\infty$.

Define the subgroup $H^*$ to be the group generated in $G(P)$ by $H$
and $\stab_H(\infty)^*$.  It follows exactly as in \cite{CL} that the
group $H^*$ is geometrically finite and leaves invariant a convex set
$C^*$ which is slightly larger than the $H^*$-orbit of $C^+ \cup
N_\infty$; moreover the set $X$ embeds into $C^*/H^*$. By construction
$H^*$ has a full rank cusp.\qed 

\subsection{}

We now discuss the existence of all right polyhedra in hyperbolic spaces. 
It is well-known that such polyhedra cannot exist for large dimensions.

We first fix some notation.  Let $\cal P$ be a convex polyhedron in ${\bf H}^n$
with a finite number of co-dimension one faces, 
all of whose dihedral angles are integer submultiples of $\pi$. 
Denote by $G({\cal P})$ (resp.
$G^+(P)$) the group generated by
reflections in the codimension one faces of $\cal P$ 
(resp. subgroup of index $2$
consisting of orientation-preserving isometries). $G(P)$ is discrete
by Poincar\'e's Theorem (see \cite{R} Chapter 7 or \cite{VS}).\\[\baselineskip]
{\leftline{\bf Ideal 24-cell in hyperbolic 4-space:}
\\[\baselineskip]
Let $P$ denote the all right ideal 24 cell $P$ in ${\bf H}^4$
(cf. \cite{R}, Example 6 p. 273, p. 509).  
We have the following lemma (see
\cite{RT}) which records some arithmetic data associated to $G(P)$.

\begin{lemma}
\label{24cell}
$G^+(P)$ is an arithmetic lattice in 
$\Isom_+({\bf H}^4)$. It is
a subgroup of finite index in $\SO_0(f_4;{\bf Z})$. \qed\end{lemma}

{\leftline{\bf A compact all right 120-cell in hyperbolic 4-space:}}

\medskip

\noindent Let $D$ denote the regular $120$-cell in ${\bf H}^4$ with all right
dihedral angles, see \cite{Co}. This has as faces 3-dimensional
all right dodecahedra. $D$ is built from 14400 copies
of the $\{4,3,3,5\}$ Coxeter simplex, $\Sigma$ in ${\bf H}^4$. 
We fix the following notation
to be used throughout. $\cal O$ will denote the ring of integers in 
${\bf Q}(\sqrt{5})$, and $\tau$ will denote the non-trivial Galois 
automorphism of ${\bf Q}(\sqrt{5})$.

\begin{lemma}
\label{120cell}
The group $G^+(D)$ is an arithmetic lattice in 
$\Isom_+({\bf H}^4)$. It is
commensurable with the group $\SO_0(f;{\cal O})$ where $f$ is the 
5-dimensional quadratic form $<1,1,1,1,-\phi>$, and $\phi = 
{{1 + \sqrt{5}}\over 2}$.
\end{lemma}
\demo By Vinberg's criteria \cite{Vi}, the group generated by
reflections in the faces of $\Sigma$ is arithmetic. By the remarks
above, $G^+(D)$ is also arithmetic.  The description given follows
from \cite{Bu}, see also \cite{VS} p. 224.\qed\\[\baselineskip]

{\leftline{\bf An all right ideal polyhedron in hyperbolic 6-space:}}

\medskip

\noindent In ${\bf H}^6$ there is a 
simplex $\Sigma$ with one ideal vertex given
by the following Coxeter diagram (see \cite{R} p. 301).

\begin{center}
\begin{picture}(170,40)

\put(1,10){\line(1,0){150}}

\put(62 , 8){\line(0,1){30}}

\put(0,7){$\bullet$}
\put(30,7){$\bullet$}
\put(60,7){$\bullet$}
\put(90,7){$\bullet$}
\put(120,7){$\bullet$}
\put(150,7){$\bullet$}

\put(60,37){$\bullet$}

\put(135,12){$4$}

\end{picture}
\end{center}
\begin{center}
{\bf Figure 1}
\end{center}
Notice that deleting the right most vertex of this Coxeter symbol 
gives an irreducible diagram for a finite Coxeter group, namely $E_6$. This
group has order $2^7.3^4.5$.

We will make use of the following.

\begin{lemma}
\label{poly6dim}
(i)~$G^+(\Sigma) = \SO_0(f_6;{\bf Z})$.  \\[\baselineskip] (ii)~There
is an all right polyhedron $Q$ built from $2^7.3^4.5$ copies of
$\Sigma$. In particular the reflection group $G(Q)$ is commensurable
with $\SO_0(f_6;{\bf Z})$.\end{lemma} \demo The first part is due to
Vinberg \cite{Vi2}, and also discussed in \cite{R} p. 301.  For the
second part, as noted above, if one deletes the face $F$ of the
hyperbolic simplex corresponding to the right hand vertex to the given
Coxeter diagram, the remaining reflection planes pass through a single
(finite) vertex and these reflections generate the finite Coxeter
group $E_6$. Take all the translates of the simplex by this group;
this yields a polyhedron whose faces all correspond to copies of
$F$. Two such copies meet at an angle which is twice the angle of the
reflection plane of the hyperbolic simplex which lies between them.
One sees from the Coxeter diagram that the plane $F$ makes angles
$\pi/2$ and $\pi/4$ with the other faces of the hyperbolic simplex, so
the resulting polyhedron is all right as required.\qed\\[\baselineskip] 
{\bf Remark.} The polyhedron is finite covolume
since there is only one infinite vertex: deleting the plane
corresponding to the left hand vertex of the Coxeter group is the only
way of obtaining an infinite group and this group is a $5$ dimensional
Euclidean Coxeter group. (See Theorem 7.3.1 Condition $(2)$ of
\cite{R}) The other Coxeter diagrams of \cite{R} shows that there are
ideal all right polyhedra in ${\bf H}^k$ at least for $2\leq k \leq 8$.

\section{Proof of Theorem \ref{main}}
\label{mainproof}

The subsections that follow collect the necessary material to be
used in the proof. 

\subsection{}

We need the the following standard facts.

\begin{lemma}
\label{sepsubgroup}
Let $G$ be a group and let $H < K < G$. If $G$ is $H$-subgroup separable then
$K$ is $H$-subgroup separable.
\end{lemma}
\demo 
Let $k \in K\setminus H$. Since $G$ is $H$-subgroup separable there is a 
finite index subgroup $G_0 < G$ with $H < G_0$ but $k \notin G_0$. Then
$G_0 \cap K$ is a subgroup of finite index in $K$ separating $H$ from
$k$ as required.\qed

\begin{lemma}
\label{geomfinite}
Let $G$ be a finite co-volume Kleinian group and let $K$ be a subgroup
of finite index.  If $K$ is $H$-subgroup separable for all geometrically
finite subgroups $H < K$, then $G$ is $H$-subgroup separable for all
geometrically finite subgroups $H < G$.
\end{lemma}
\demo It is a standard fact (see for example Lemma 1.1 of \cite{Sc})
that if $K$ is subgroup separable and $K \leq G$ with $K$ being of
finite index, then $G$ is subgroup separable. The proof of this result
applies verbatim to Lemma \ref{geomfinite} after one notes that the
property of being geometrically finite is preserved by super- and sub-
groups of finite index. \qed

\begin{lemma}
\label{geomfinite2}
Let $\Gamma$ be a discrete subgroup of $\Isom_+({\bf H}^n)$ of finite
co-volume and $\Gamma_0$ a geometrically finite
subgroup of $\Gamma$ fixing a totally geodesic copy of hyperbolic 3-space 
in ${\bf H}^3$.
Then $\Gamma_0$ is geometrically finite as a subgroup of $\Isom_+({\bf H}^n)$.
\end{lemma}
\demo Given a geometrically finite hyperbolic (n-1)-orbifold 
${\bf H}^{n-1}/G$
this can be seen to be geometrically finite as a quotient of
${\bf H}^n$ by observing that an $\epsilon$-neighbourhood of the 
core in ${\bf H}^n$ is isometric
to  (core in ${\bf H}^{n-1}) \times I$ since $G$ fixes a genuine
co-dimension one geodesic sub-hyperbolic space.
That the $n$-dimensional core has finite volume now follows from the
fact that the $(n-1)$-dimensional core does. The proof of the
statement we require follows from this and induction. Alternatively
one could use Lemma \ref{GF2} and observe that the properties
of conical limit points
and bounded parabolic fixed points will be preserved.\qed

\subsection{}

The key lemma is the following. 

\begin{lemma}
\label{subgpbianchi}
Let $f$ be the quadratic form $<1,1,1,1,1,1,-1>$. Then
for all $d$,  $\SO(f;{\bf Z})$ contains a group 
$G_d$ which is conjugate to a subgroup of finite index in the Bianchi
group $\PSL(2,O_d)$.
\end{lemma}

The proof requires an additional lemma.  Assume that $j$ is a diagonal
quaternary quadratic form with integer coefficients of signature
$(3,1)$; so that $j$ is equivalent over $\bf R$ to the form
$<1,1,1,-1>$. Let $a \in {\bf Z}$ be a square-free positive integer
and consider the seven dimensional form $j_a = <a,a,a> \oplus ~j$,
where $\oplus$ denotes orthogonal sum. Being more precise, if we
consider the 7-dimensional $\bf Q$-vector space $V$ equipped with the
form $j_a$ there is a natural $4$-dimensional subspace $V_0$ for which
the restriction of the form is $j$. Using this it easily follows that,

\begin{lemma}
\label{subgroup}
In the notation above, the group $\SO(j;{\bf Z})$ is a subgroup of 
$\SO(j_a;{\bf Z})$.\qed\end{lemma}

\leftline{{\bf Proof of Lemma \ref{subgpbianchi}}}

\noindent Let $p_d$ be the quaternary form $<d,1,1,-1>$. Notice that
this form represents $0$ non-trivially, and hence the corresponding
arithmetic group $\SO_0(p_d;{\bf Z})$ is non-cocompact. By Theorem
\ref{macr} for example, this implies that $\SO_0(p_d;{\bf Z})$ is
commensurable with some conjugate of an appropriate image of the
Bianchi group $\PSL(2,O_d)$. The key claim is that $q_d = <d,d,d>
\oplus~ p_d$ is equivalent over $\bf Q$ to the form $f$.

Assuming this claim for the moment, by Lemma \ref{equivalent} we deduce that
there exists $R_d \in \GL(7,{\bf Q})$ such that 
$R_d\SO(q_d;{\bf Z})R_d^{-1}$ and $\SO(f;{\bf Z})$ are commensurable. This
together with Lemma \ref{subgroup} gives the required group $G_d$.

To prove the claim, since every positive integer can be written as
the sum of four squares, write $d = w^2 + x^2 + y^2 + z^2$.
Let $A_d$ be the $7 \times 7$ matrix 
$$\pmatrix{
w & x & y & z & 0 & 0 & 0 \cr -x & w & -z & y & 0 & 0 & 0 \cr
-y & z & w & -x & 0 & 0 & 0 \cr -z & -y & x & w & 0 & 0 & 0 \cr
0 & 0 & 0 & 0 & 1 & 0 & 0 \cr 0 & 0 & 0 & 0 & 0 & 1 & 0 \cr
0 & 0 & 0 & 0 & 0 & 0 & 1\cr}
$$
Note $A_d$ has determinant $d^2$, so is in $\GL(7,{\bf Q})$.
Let $F$ the diagonal matrix associated to the form $f$ and $Q_d$ be
the $7\times 7$ diagonal matrix of the form $<d,d,d,d,1,1,-1>$ (i.e. of
$<d,d,d> \oplus ~p_d$).  Then a direct check shows that $A_dFA_d^t = Q_d$
as is required.\qed\\[\baselineskip] 
{\bf Remark:}~The appearance of the matrix $A_d$ is described by the following.
Let $\Omega = \{w + xi + yj + zij ~:~w,x,y,z \in {\bf Z}\}$ denote the ring
of integral Hamiltonian quaternions. If $\alpha = w + xi + yj 
+ zij \in \Omega$,
then the norm of $\alpha$ is $d = w^2 + x^2 + y^2 + z^2$. By considering
the representation of $\Omega$ into $M(4,{\bf Q})$ determined by
the right action of $\Omega$ on itself, $\alpha$ is mapped to the $4 \times 4$
block of the matrix $A_d$, and $\overline{\alpha} = w - xi - yj - zij$,
to the transpose of this block.\\[\baselineskip]
{\bf Proof of Theorem \ref{main}.}  
By Theorem \ref{allright} and Lemma \ref{poly6dim} we
deduce that $\SO_0(f;{\bf Z})$ is $H$-subgroup separable on all its
geometrically finite subgroups $H$.  By Lemma \ref{subgpbianchi} the
groups $G_d$ are subgroups of $\SO_0(f;{\bf Z})$, and so by Lemma
\ref{geomfinite2} they, and all their geometrically finite subgroups
(as groups acting on ${\bf H}^3$) are geometrically finite subgroups
of $\SO_0(f;{\bf Z})$ (acting on ${\bf H}^6$). Hence Lemma
\ref{sepsubgroup} shows that $G_d$ is $H$-subgroup separable for all
geometrically finite subgroups $H$ of $G_d$.  Lemma \ref{geomfinite}
allows us to promote this subgroup separability to the groups
$\PSL(2,O_d)$.  This proves Theorem \ref{main}.\qed

\subsection{The cocompact case.}

We can extend the techniques used in the proof of Theorem 
\ref{main} to cocompact groups. The crucial lemma is:

\begin{lemma}
\label{infinitelymany}
Let $f$ be the quadratic form $<1,1,1,1,-1>$, and $q$ 
the quadratic form $<1,1,1,1,-\phi>$ (as in Lemma \ref{120cell}). Then,
\smallskip
\item{1.}~$\SO_0(f;{\bf Z})$
contains infinitely many commensurability classes of cocompact arithmetic
Kleinian groups. Furthermore these can be chosen to be incommensurable
with any group generated by reflections in dimension $3$.
\smallskip
\item{2.}~$\SO_0(q;{\cal O})$
contains infinitely many commensurability classes of cocompact arithmetic
Kleinian groups. Furthermore these can be chosen to be incommensurable
with any group generated by reflections in dimension $3$.
\end{lemma}

The lemma will be proved in \S \ref{firstpartinfinitelymany}, assuming
it we establish Theorem \ref{cocompact}.\\[\baselineskip]
{\bf Proof of Theorem \ref{cocompact}.} 
The proof is identical to that of the implication of \ref{main} from
Lemma \ref{subgpbianchi}. In this case we use
Theorem \ref{allright}, Lemma \ref{24cell} and Lemma \ref{120cell}
to get the appropriate all right reflection group, and commensurability
with the special orthogonal groups in question.\qed

\section{Preliminaries for Lemma \ref{infinitelymany}}

We use this section to record facts about equivalence of quadratic forms 
that we will need, (see \cite{La}).

Let $K$ denote either a number field or a completion of
a number field, and
$q$ a non-singular quadratic form defined over $K$ with associated
symmetric matrix $Q$. The {\em determinant} of $q$ is the element 
$d(q) = \det (Q)\dot{K}^{2}$, where $\dot{K}$ are the invertible 
elements in $K$. It is not hard to see that $d(q)$ is
an invariant of the equivalence class of $q$.

The {\em Hasse invariant} (see \cite{La}, p. 122) of a non-singular
diagonal form $<a_1,a_2, \ldots ,a_n>$ with coefficients in $K$
 is an element in the Brauer group $B(K)$, namely
$$s(q) = \prod_{i < j} \biggl({{a_i,a_j} \over K}\biggr)$$
where $\biggl({{a_i,a_j} \over K}\biggr)$ describes a quaternion
algebra over $K$, and the multiplication is that in $B(K)$,
see \cite{La}, Chapter 4.

Every non-singular form over $K$ is equivalent over $K$ to a diagonal
one, and the definition of the Hasse invariant is extended to non-diagonal
forms by simply defining it to be the Hasse invariant of a diagonalization
(that this is well-defined follows from \cite{La}, p. 122)
The following theorem is important to us. It is called the
``Weak Hasse-Minkowski Principle'' in \cite{La}, p. 168.  We state it
in the case when $K$ is a number field.

\begin{theorem}
\label{weakHM}
Let $q_1$ and $q_2$ be non-singular quadratic forms of the same dimension,
defined over $K$ with the property that if $\sigma$ is a real
embedding of $K$ the forms $q_1^\sigma$ and $q_2^\sigma$ have the
same signature over $\bf R$. Then $q_1$ is
equivalent to $q_2$ over $K$ if and only if $d(q_1) = d(q_2)$ and 
$s(q_1) = s(q_2)$ over all non-archimedean completions of $K$.\qed\end{theorem}

Note that if $d(q_1)$ = $d(q_2)$ (resp. $s(q_1) = s(q_2)$) then
the same holds locally.

\section{Proof of Lemma \ref{infinitelymany}(1)}
\label{firstpartinfinitelymany}

In this section we give the proof of the first part of Lemma 
\ref{infinitelymany}. The method of proof of the third part is the same,
however some additional algebraic complexities are involved
since we are working over the field
${\bf Q}(\sqrt{5})$. This is dealt with in the next section.

\subsection{}

Assume that $j$ is a diagonal quaternary quadratic form with integer
coefficients of signature $(3,1)$; so that $j$ is equivalent over $\bf
R$ to the form $<1,1,1,-1>$. Let $a \in {\bf Z}$ be a square-free positive
integer and consider the five dimensional form $j_a = <a> \oplus ~j$, where
$\oplus$ denotes orthogonal sum. As in \S 4, if we consider the
5-dimensional $\bf Q$-vector space $V$ equipped with the form $j_a$
there is a natural $4$-dimensional subspace $V_0$ for which the
restriction of the form is $j$. As before it easily follows that,

\begin{lemma}
\label{subgroup2}
In the notation above, the group $\SO(j;{\bf Z})$ is a subgroup of 
$\SO(j_a;{\bf Z})$.\qed\end{lemma}

\noindent We begin the proof of the first claim in
Lemma \ref{infinitelymany}. Let $p_d$ denote the quadratic
form $<1,1,1,-d>$. This has signature $(3,1)$, and as discussed in \S \ref{preliminaries}
gives arithmetic Kleinian groups. We have the following classical result
from number theory, see \cite{La}, pp. 173--174.

\begin{theorem}
\label{3squares}
Let $d$ be a positive integer. Then $d$ is the sum of three squares if and
only if $d$ is not of the form $4^t(8k - 1)$.\qed\end{theorem}

\noindent Choose a square-free positive integer 
$d  = -1~ \mod~ 8$, and let $p_d$ the form $<1,1,1,-d>$. Since
a non-trivial rational solution $p_d(x) = 0$ can be easily promoted to
an integral solution, 
Theorem \ref{3squares}, shows this form
does not represent $0$ non-trivially over $\bf Q$. Hence 
the arithmetic Kleinian groups $\SO_0(p_d;{\bf
Z})$ are cocompact.  By Theorem \ref{macr},
to get the Kleinian groups to be incommensurable, 
we simply insist further that $d$ is a prime. By Dirichlet's
Theorem there are infinitely many such primes. With these remarks,
it follows from Theorem \ref{macr} that the groups $\SO(p_d;{\bf Z})$ are all
incommensurable. The first part of Lemma \ref{infinitelymany} will follow
from.

\begin{lemma}
\label{hasse_calcs}
Let $q_d = <d> \oplus ~ p_d$. Then $q_d$ is equivalent over $\bf Q$ to $f$.
\end{lemma}
\demo  The two forms are $5$-dimensional, and it is easy to see that the 
forms have signature $(4,1)$ over $\bf R$. Further since the determinants
are $-1\dot{\bf Q}^{2}$, they will have the same local determinants.
We shall show that the forms have the same Hasse invariants over $\bf Q$
from which it follows they have the same local Hasse invariants. Theorem
\ref{weakHM} completes the proof.

Consider the form $f$ first of all.  It is easy to see that all the
contributing terms to the product are either 
$\biggl({{1,1} \over {\bf Q}}\biggr)$ or 
$\biggl({{-1,1} \over {\bf Q}}\biggr)$. Both of these are isomorphic to
the quaternion algebra of $2 \times 2$ matrices over $\bf Q$, see \cite{La}
Chapter 3. These represent the trivial element in the Brauer group of
$\bf Q$, and so $s(f)$ is trivial.

For $q_d$ the contributing terms are 
$$\biggl({{1,1} \over {\bf Q}}\biggr),~\biggl({{1,d} \over {\bf Q}}\biggr),~
\biggl({{1,-d} \over {\bf Q}}\biggr), \biggl({{d,-d} \over {\bf Q}}\biggr).$$
From \cite{La} Chapter 3, in particular p. 60, it follows that all these
quaternion algebras are again isomorphic to
the quaternion algebra of $2 \times 2$ matrices over $\bf Q$, and so as above
$s(q_d)$ is trivial. This completes the proof.\qed\\[\baselineskip]
{\bf Remark:}~The proof of this Lemma can also be done directly as in the
case of Lemma \ref{subgpbianchi}. We include this version, as it may be useful
as a guide to the proof of Lemma \ref{hasse_calcs2}.\\[\baselineskip]
To complete the proof of the first claim in Lemma \ref{infinitelymany}
we proceed as follows---entirely analogous to the argument in the
proof of \ref{main}.  By Lemma \ref{subgroup2}, $\SO(p_d;{\bf Z})$
is a subgroup of $\SO(q_d;{\bf Z})$. By Lemma \ref{hasse_calcs} and
Lemma \ref{equivalent} we can conjugate to obtain a group $G_d < SO(f;{\bf Z})$
which is conjugate to a subgroup of finite index in $\SO(p_d;{\bf Z})$.
Finally, the groups
constructed are not commensurable with groups generated by reflections
for $d$ large enough. This follows from work of Nikulin, \cite{Ni}.

Briefly it is shown in \cite{Ni} that if we fix the field of definition
for a reflection group (in this case $\bf Q$) there are
only finitely many commensurability classes of arithmetic Kleinian groups
commensurable with a group generated by reflections.
This completes the proof.\qed\\[\baselineskip]
{\bf Remarks.}\\
1.~The argument given in Lemma \ref{infinitelymany} also applies in the
case of the Bianchi groups.\\
2.~An example of an explicit ``well known'' co-compact group that is covered 
by our techniques arises in the choice of $d = 7$.
It follows from \cite{MR} that the arithmetic Kleinian
group $\Gamma$ arising as an index 2 subgroup in the group generated by 
reflections in the faces of the tetrahedron $T_6[2,3,4;2,3,4]$
is commensurable with the group $\SO(p_7;{\bf Z})$.

Also commensurable with $\SO(p_7;{\bf Z})$ is the fundamental group of
a certain non-Haken hyperbolic 3-manifold obtained by filling on a
once-punctured torus bundle. Briefly, let $M$ denote the once-punctured
torus bundle whose monodromy is $R^2L^2$ (in terms of the usual
$RL$-factorization, see \cite{CJR} for instance). As is well-known,
$M$ contains
no closed embedded essential surfaces. Thus a Haken manifold
can be created by Dehn filling on $M$, only by filling along a boundary slope
of $M$.
Fixing a framing for
the boundary torus, the manifold $M_0$ obtained by $4/1$-Dehn filling
on $M$ is hyperbolic, of volume approximately $2.666744783449061\ldots$,
and non-Haken (the boundary slopes can be deduced
from \cite{CJR}). A calculation with Snap (see \cite{CGHN} for a
discussion of this program) shows $M_0$ is arithmetic with the same
invariant data as the group $\Gamma$ above. Hence $\pi_1(M_0)$ is
separable on geometrically finite subgroups.

\section{Proof of Lemma \ref{infinitelymany}(2)}
\label{proofofinfinitelymany}

To handle the second part, we proceed in a similar way to 
\S \ref{firstpartinfinitelymany}.

Thus, assume that $j$ is a diagonal quaternary quadratic form with 
coefficients in ${\cal O}$, of signature $(3,1)$ at the
identity and signature $(4,0)$ on applying $\tau$.
Let $a \in {\cal O}$ be postive at both the identity embedding and $\tau$.
Consider the five dimensional form $j_a = <a> \oplus ~j$, where
$\oplus$ denotes orthogonal sum. As above, we mean
5-dimensional ${\bf Q}(\sqrt{5})$-vector space $V$ equipped with the form $j_a$
there is a natural $4$-dimensional subspace $V_0$ for which the
restriction of the form is $j$. 
We have
the following consequence of the discussion in \S 2.2. 
\begin{lemma}
\label{sqrt5}
Let $a \in {\cal O}$ have the property that it is
square-free (as an element of ${\cal O}$), $a < 0$
but $\tau (a) > 0$. Define the $n+1$-dimensional form 
$$f_{n,a} = <1,1,\ldots,1, a>.$$
Then $\SO_0(f_{n,a};{\cal O})$ defines a cocompact arithmetic
subgroup of $\Isom_+({\bf H}^n)$.\qed\end{lemma}

With this, we deduce as in \S 6.1, 

\begin{lemma}
\label{subgroup3}
In the notation above, the group $\SO(j;{\cal O})$ is a cocompact
subgroup of $\SO(j_a;{\cal O})$.\qed\end{lemma}

\subsection{}

We need to recall some basic number theory in $\cal O$.
$\cal O$ is a principal ideal domain, so every ideal of $\cal O$
has the form $<t>$ for some $t\in \cal O$.

\begin{lemma}
\label{generator}
Let $I \subset \cal O$ be a non-trivial ideal. Then $I$ can be
generated by an element $t$ where $t < 0$ and $\tau (t) >
0$.\end{lemma}
\demo Assume $I = <x>$. The argument is really a consequence of the existence
of units of all possible signatures in $\cal O$. For example if $x$
is positive at both the identity and $\tau$, then 
$t = x({{1 + \sqrt{5}}\over 2})$ satisfies the requirements. As 
${{1 + \sqrt{5}}\over 2}$ is a unit it does not change the ideal.
\qed\\[\baselineskip]
Define ${\bf P}$ to be the set of primes in ${\bf Q}(\sqrt{5})$ 
lying over the set of rational primes $\cal Q$ satisfying:
\begin{enumerate}
\item if $q \in {\cal Q}$ then $q$ is unramified in 
${\bf Q}(\sqrt{{1+\sqrt{5}}\over 2})$,
\item the splitting type of $q$ in ${\bf Q}(\sqrt{{1+\sqrt{5}}\over 2})$
contains a prime of degree 1, ie. there is a 
${\bf Q}(\sqrt{{1+\sqrt{5}}\over 2})$-prime $Q$ with $Q|q$ and the norm
of $Q$ is $q$.
\end{enumerate}

Note by the Tchebotarev density theorem (see \cite{Ln} for example), the
set $\bf P$ is infinite, since the set of rational
primes that split completely in 
${\bf Q}(\sqrt{{1+\sqrt{5}}\over 2})$ is infinite.
Fix $\phi = {{1+\sqrt{5}}\over 2}$ in what follows.

\begin{lemma}
\label{square}
Let ${\cal P} \in {\bf P}$. Then under the canonical reduction map 
${\cal O} \rightarrow {\cal O}/{\cal P}$, $\phi$ 
is a square.\end{lemma}
\demo  
Let $\cal R$ denote the ring of integers in
${\bf Q}(\sqrt{{1+\sqrt{5}}\over 2})$. If $p \in {\cal Q}$, then since 
${\bf Q}(\sqrt{5})$ is galois, $p$
splits completely in ${\bf Q}(\sqrt{5})$. By assumption there is
a degree one prime $P \subset {\cal R}$ with $P | p{\cal R}$.
Let ${\cal P} \in {\bf P}$ be a ${\bf Q}(\sqrt{5})$-prime divisible by $P$.
Then the ramification theory of primes in extensions (see \cite{Ln}), together
with the above properties give,
$${\cal R}/P \cong  {\cal O}/{\cal P} \cong {\bf F}_p.$$
Hence by definition, $\phi$ has a square-root
upon reduction $\mod ~{\cal P}$ as required.
\qed\\[\baselineskip]

Define the following collection of quadratic forms over ${\bf Q}(\sqrt{5})$.
Let ${\cal P} \in {\bf P}$ be generated by an element $\pi$ satisfying
the conclusion of Lemma \ref{generator}, that is 
$$\pi < 0 ~\hbox{and}~ \tau (\pi) > 0.$$
Define
$p_\pi = <1,1,1,\pi>$. 
Thus, by Lemma \ref{sqrt5}, $\SO(p_\pi;{\cal O})$ determines
a cocompact arithmetic Kleinian group. Furthermore this infinite
collection of groups
are all mutually incommensurable. This follows
\cite{B} or \cite{MR} (this is the generalization of Theorem \ref{macr}
to ${\bf Q}(\sqrt{5})$).

Define the $5$-dimensional form
$q_\pi = <-\pi\phi> \oplus ~p_\pi$. 
Since $\pi$ (resp. $\tau (\pi)$) is negative (resp. positive)
and $\phi$ (resp. $\tau (\phi)$) is positive (resp. negative), 
both $-\pi\phi$ and $-\tau (\pi)\tau (\phi)$ are positive.
Thus $q_\pi$ has signature $(4,1)$ at the identity and $(5,0)$ at $\tau$.
Therefore, by Lemma \ref{subgroup3},
$\SO(q_\pi;{\cal O})$ determines a cocompact arithmetic subgroup
of $\Isom_+({\bf H}^4)$, which contains $\SO(p_\pi;{\cal O})$.

In what follows we quote freely from the theory of quaternion algebras,
see \cite{La} Chapter 3 and \cite{Vig}.

We will need the following.

\begin{lemma}
\label{hilbert}
Let $k = {\bf Q}(\sqrt{5})$. The quaternion algebra $B =
\biggl({{\phi,\pi} \over k}\biggr)$ is isomorphic to
$M(2,k)$.\end{lemma} \demo Note that by choice of $\pi$, $B$ is
unramified at both the identity embedding and $\tau$. Now Lemma
\ref{square} implies $\phi$ is a square upon reduction $\mod ~\pi$.
It follows that the norm form of $B$, namely $<1,-\phi,-\pi,\phi\pi>$
is isotropic over $k_{<\pi>}$, since it is isotropic over the residue
class field (see \cite{La} Chapter 6). The only primes that can
ramify $B$ are $<\pi>$ and the unique prime above $2$ in $k$ (cf. \cite{Vig}
or \cite{La} Chapter 6). Furthermore since the cardinality of
the ramification set is even (\cite{Vig}), we deduce that $B$
is unramified at the prime above $2$. Hence $B$ is unramified
everywhere locally, and so is isomorphic to $M(2,k)$.\qed

\begin{lemma}
\label{hasse_calcs2}
In the notation above, $q_\pi$ is equivalent over ${\bf Q}(\sqrt{5})$ to $q$.
\end{lemma}

\demo  Let $k = {\bf Q}(\sqrt{5})$.
The two forms are $5$-dimensional, and as noted both
forms have signature $(4,1)$ at the identity embedding of $k$, and signature
$(5,0)$ at $\tau$. Further since the determinants
are $-\phi\dot {k}^{2}$,
they will have the same local determinants.
We shall show that the forms have the same Hasse invariants over $k$
from which it follows they have the same local Hasse invariants. Theorem
\ref{weakHM} completes the proof.

Consider the form $q$ first of all.  It is easy to see that all the
contributing terms to the product are either 
$\biggl({{1,1} \over k}\biggr)$ or $\biggl({{1,-\phi} \over k}\biggr)$.
Both of these are isomorphic to
the quaternion algebra of $2 \times 2$ matrices over $k$, see \cite{La}
Chapter 3 (in particular p. 60). 
These represent the trivial element in the Brauer group of
$k$, and so $s(f)$ is trivial.

For $q_\pi$ the contributing terms are 
$$\biggl({{1,1} \over k}\biggr),~\biggl({{1,-\phi\pi} \over k}\biggr),~
\biggl({{1,\pi} \over k}\biggr), \biggl({{\pi,-\phi\pi} \over k}\biggr).$$
As above, it follows from \cite{La} Chapter 3, p. 60, that all but the
last Hilbert symbol represent 
quaternion algebras isomorphic to
the quaternion algebra of $2 \times 2$ matrices over $k$.

Standard Hilbert symbol manipulations (\cite{La} Chapter 3) 
imply this last algebra is isomorphic
to one with Hilbert Symbol $\biggl({{\phi,\pi} \over k}\biggr)$. But
Lemma \ref{hilbert} implies this quaternion algebra is the matrix algebra 
again. Hence $s(q_\pi)$ is also trivial. 
This completes the proof.\qed\\[\baselineskip]
To complete the proof of Lemma \ref{infinitelymany}(3)
we proceed as follows.  As noted all groups are cocompact, and
Lemma \ref{subgroup3} implies
that $\SO(p_\pi;{\cal O})$
is a subgroup of $\SO(q_\pi;{\cal O})$. By Lemma \ref{hasse_calcs} and
Lemma \ref{equivalent} we can conjugate to obtain a group 
$G_\pi < \SO(f;{\cal O})$
which is conjugate to a subgroup of finite index in $\SO(p_\pi;{\cal O})$,
The final part is to deduce that
infinitely many of these are not commensurable with groups generated
by reflections, and again this follows from \cite{Ni}, with the
ground field in this case being ${\bf Q}(\sqrt{5})$.
\qed\\[\baselineskip]

Corollary \ref{cox} is deduced 
by choosing $\pi = (3 - 2\sqrt{5})$ in the construction above.
Note $<\pi>$ is a prime in ${\bf Q}(\sqrt{5})$ of norm $11$, and as is
easily checked $11$ is unramified in ${\bf Q}(\sqrt{\phi})$ and 
also satisfies the second condition in the definition of $\bf P$ in \S 4.2.
It follows from \cite{MR} that the arithmetic Kleinian
group $\Gamma_0$ is commensurable with $\SO(p_\pi;{\cal O})$.
\qed\\[\baselineskip]

Corollary \ref{SWD} is deduced in a similar manner. The Seifert-Weber
dodecahedral space is constructed from $120$ copies of the tetrahedron
$T_4[2,2,5;2,3,5]$. If we let $\Gamma$ denote the
group generated by reflections in faces of this tetrahedron, the results of
\cite{MR} imply $\Gamma$ is commensurable (up to conjugacy) with the group
$\SO(f;{\cal O})$ where $f$ is the form $<1,1,1,-1-2\sqrt{5}>$.
Now $<-1 - 2\sqrt{5}>$ generates a prime ideal of norm $19$, and,  as is 
easily checked, lies in $\bf P$. The result now follows.
\qed

\section{Application.}

Theorem \ref{main} allows to address the following question of
Lubotzky raised in \cite{Lu}. We require some terminology.
The profinite topology on a group $G$ is defined by proclaiming all
finite index subgroups of $G$ to be a basis of open neighbourhoods of
the identity. Let $\hat{G}$ denote the profinite completion of $G$.

An obvious subgroup of $\PSL(2,O_d)$ is $\PSL(2,{\bf Z})$, and in the
context of the Congruence Kernel, Lubotzky \cite{Lu} asked the
following question: \\[\baselineskip] 
{\bf Question}: Is the induced map:
 $$\eta~:~\widehat{\PSL}(2,{\bf Z}) \rightarrow
\widehat{\PSL}(2,O_d)$$ 
injective?   \\[\baselineskip]

Since open subgroups are closed in the profinite topology, it is not
hard to see that $G$ is $H$-subgroup separable if and only if $H$ is
closed in the profinite topology on $G$. Thus an equivalent formulation of
Theorem \ref{main} is that  if $H$ is a geometrically finite subgroup 
of $\PSL(2,O_d)$, then on passing to the profinite completion 
$\widehat{\PSL}(2,O_d)$,  we require that the only points of $\PSL(2,O_d)$ 
in the closure of $H$ (in $\widehat{\PSL}(2,O_d)$)  are points of $H$.
This formulation can be can be used to give an affirmative answer.

\begin{theorem}
The map $\eta_d$ is injective for all $\PSL(2,O_d)$
\end{theorem}
\demo The group $\PSL(2,O_d)$ is residually finite, so that it (and hence
$\PSL(2,{\bf Z}) )$  embeds into the  profinite completion $\widehat{\PSL}(2,O_d)$.
Thus we may take the closure of $\PSL(2,{\bf Z}) \subset \widehat{\PSL}(2,O_d)$
to form a completion, denoted $\overline{\PSL}(2,{\bf Z})$. This completion is
obviously embedded in  $\widehat{\PSL}(2,O_d)$, but it is potentially coarser
than the genuine profinite completion  $\widehat{\PSL}(2,{\bf Z})$.

However any finitely generated subgroup $H$ (in particular any subgroup of finite index)
in  $\PSL(2,{\bf Z})$ is geometrically finite and thus by Theorem \ref{main}
separable in $\PSL(2,O_d)$
and it follows easily that one can find a subgroup $H^*$ of finite index in
 $\PSL(2,O_d)$ with the property that $H^* \cap   \PSL(2,{\bf Z}) = H$ so that in
fact the map   $\widehat{\PSL}(2,{\bf Z}) \rightarrow \overline{\PSL}(2,{\bf Z})$ is 
a homeomorphism, proving the theorem. $\Box$

\noindent Department of Mathematics,\\
University of Melbourne\\
Parkville, VIC 3010 Australia\\[\baselineskip] 
Department of Mathematics,\\
University of California\\
Santa Barbara, CA 93106\\[\baselineskip] 
Department of Mathematics,\\
University of Texas\\
Austin, TX 78712\\

\end{document}